\documentclass[10pt]{elsarticle}
\journal{}
\usepackage{amsmath}
\usepackage{indentfirst}
\usepackage{amsfonts,amssymb}
\usepackage{mathrsfs}
\usepackage{amsthm}
\usepackage{lineno,hyperref}
\modulolinenumbers[5]
\newtheorem{Definition}{Definition}[section]
\newtheorem{Theorem}{Theorem}[section]
\newtheorem{Lemma}{Lemma}[section]

\numberwithin{equation}{section}

\begin{document}
	
	\begin{frontmatter}
		
		\title{Dynamics of an imprecise stochastic Lotka–Volterra food chain chemostat model with L\'{e}vy jumps}

		\author{Fei Sun\corref{mycorrespondingauthor}}
		\cortext[mycorrespondingauthor]{Corresponding author}
		\address{School of Mathematics and Computational Science, Wuyi University, Jiangmen 529020, China}
		\ead{sunfei@whu.edu.cn (fsun.sci@outlook.com)}

\begin{abstract}Population dynamics are often affected by sudden environmental perturbations. Parameters of stochastic models are often imprecise due to various uncertainties. In this paper, we formulate a stochastic Lotka–Volterra food chain chemostat model that includes L\'{e}vy jumps and interval parameters. Firstly, we prove the existence and uniqueness of the positive solution. Moreover, the threshold between extinction and persistence in the mean of the microorganism is obtained. Finally, some simulations are carried out to demonstrate our theoretical results.
\end{abstract}


\end{frontmatter}

\section{Introduction}
\label{sec:1}

In recent decades, the study of chemostat has become one of the main topics in mathematical biology and microbial ecology. For a better review of mathematical models on the theory of chemostat, see Smith and Waltman \cite{1} and Chen et al. \cite{2}. 

Considering the chemostat is inevitably affected by environmental noise.
May \cite{11} pointed out that some parameters involved in the system should exhibit random fluctuation to a greater or lesser extent. In order to capture essential feature of the continuous culture process of microorganisms in the chemostat, some researchers have studied the dynamics of chemostat models driven by white noise (see e.g. \cite{12}-\cite{19} as well as there references).

Most population systems assume that model parameters are accurately known. However, the sudden environmental perturbations may bring substantial social and economic losses. For example, the recent COVID-19 has a serious impact on the world. It is more realistic to study the population dynamics with imprecise parameters. Panja et al. \cite{23} studied a cholera epidemic model with imprecise numbers and discussed the stability condition of equilibrium points of the system. Das and Pal \cite{24} analyzed the stability of the system and solved the optimal control problem by introducing an imprecise SIR model. Other studies on imprecise parameters include those of \cite{20}-\cite{222}, and the references therein.

The main focus of this paper is Dynamics of an imprecise stochastic Lotka–Volterra food chain chemostat model with L\'{e}vy jumps. To this end, we first introduce the imprecise stochastic Lotka–Volterra food chain chemostat model. With the help of Lyapunov functions, we prove the existence and uniqueness of the positive solution. Further, the threshold between extinction and persistence in the mean of the microorganism is obtained.

The remainder of this paper is organized as follows. In Sect.~\ref{sec:2}, we introduce the basic models. In Sect.~\ref{sec:3}, the unique global positive solution of the system is proved. 
The threshold between extinction and persistence in the mean of the microorganism are derived in
Sect.~\ref{sec:4} and Sect.~\ref{sec:5}. 

\section{Imprecise stochastic Lotka–Volterra food chain chemostat model}
\label{sec:2}

In this section, we introduce the imprecise stochastic Lotka–Volterra food chain chemostat model.
Let $ S(t) $, $ x(t) $ and $ y(t) $ denote the concentrations of nutrient, prey and predator at time $ t $, respectively. Then a stochastic Lotka–Volterra food chain chemostat model takes the following form \cite{27}.
\begin{equation}\label{2.1}
\left\{
\begin{array}{lcl}
dS(t) = [D(S^{0}-S(t))-\dfrac{m_{1}S(t)x(t)}{\delta_{1}}]dt + \sigma_{1}S(t)dB_{1}(t) +  \int_{\mathbb{Y}}\gamma_{1}S(t^{-})\widetilde{N}(dt,du),\\
dx(t) = [m_{1}S(t)x(t)-Dx(t)-\dfrac{m_{2}x(t)y(t)}{\delta_{2}}]dt + \sigma_{2}S(t)dB_{2}(t) +  \int_{\mathbb{Y}}\gamma_{2}x(t^{-})\widetilde{N}(dt,du),  \\
dy(t) = [m_{2}x(t)y(t)-Dy(t)]dt + \sigma_{3}y(t)dB_{3}(t) +  \int_{\mathbb{Y}}\gamma_{3}y(t^{-})\widetilde{N}(dt,du).
\end{array}  
\right.
\end{equation}

where $ S^{0} $ is the original input concentration of nutrient and $ D $ is the input rate from the feed vessel to the culture vessel, as well as the washout rate from the culture vessel to the receptacle. $ \delta_{1} $ and $ \delta_{2} $ are the yield constants for prey growth on nutrient and predator growth on prey, respectively.  $ B_{i}(t) $, $ i = 1, 2, 3 $ are standard Brownian motions. $ \sigma_{i}^{2} $ represent the intensities of  $ B_{i}(t) $. $ \lambda $ is the characteristic measure of $ N $ which is defined on a finite measurable subset $ \mathbb{Y} $ of $ (0,+\infty) $ with $ \lambda(\mathbb{Y}) < \infty $. $ \gamma_{i}(u) :\mathbb{Y} \times \Omega \rightarrow \mathbb{R}  $  are the bounded and continuous functions satisfying $ \gamma_{i}(u) > -1 $, $ i = 1, 2, 3 $.
 The $ S(t^{-}) $, $ x(t^{-}) $ and $ y(t^{-}) $ are the left limits of  $ S(t) $, $ x(t) $ and $ y(t) $, respectively. $ \widetilde{N} $ denotes the compensated random measure defined by $ \widetilde{N}(dt, du) = N(dt, du) − \lambda(du)dt $. We assume $ B_{i} $ and $ N $ are independent throughout the paper and denote $ \mathbb{R}^{d}_{+} = \{x \in \mathbb{R}^{d} : x_{i} > 0 \textrm{ for all } 1 \leq i \leq d\} $, $ \overline{\mathbb{R}}^{d}_{+}
= \{x \in \mathbb{R}^{d}: x_{i} \geq 0 \textrm{ for all } 1 \leq i \leq d\} $. If $ f(t) $ is an integrable function on $  [0, \infty) $, define $\langle f \rangle_{t} =\dfrac{1}{t} \int_{0}^{t} f(s) ds$.

Before we state the imprecise stochastic Lotka–Volterra food chain chemostat model, definitions of Interval-valued function should recalled (Pal \cite{24}).

\begin{Definition}
(Interval number) An interval number $ A $ is represented by closed interval $ [a^{l}, a^{u}] $ and defined by $ A = [a^{l}, a^{u}]  = \{x|a^{l} \leq x \leq a^{u}, x \in \mathbb{R}\},$ where $ \mathbb{R} $ is the set of real numbers and $ a^{l} $, $  a^{u} $ are the lower and upper limits of the interval numbers, respectively. The interval number $ [a, a] $ represents a real number $ a $. The arithmetic operations for any two interval numbers $ A = [a^{l}, a^{u}]$ and $ B = [b^{l}, b^{u}]$ are as follows:\\
\indent Addition: $ A+B= [a^{l}, a^{u}] + [b^{l}, b^{u}] = [a^{l}+b^{l}, a^{u}+b^{u}] $.\\
\indent Subtraction:$ A-B= [a^{l}, a^{u}] - [b^{l}, b^{u}] = [a^{l}-b^{l}, a^{u}-b^{u}] $.\\
\indent Scalar multiplication: $\alpha A= \alpha [a^{l}, a^{u}] =  [\alpha a^{l}, \alpha a^{u}]$, where $ \alpha $ is a positive real number.\\
\indent  Multiplication: $ AB= [a^{l}, a^{u}]  [b^{l}, b^{u}]= [\min \{a^{l}b^{l}, a^{u}b^{l}, a^{l}b^{u}, a^{u}b^{u} \}, \max \{a^{l}b^{l}, a^{u}b^{l}, a^{l}b^{u}, a^{u}b^{u} \}  ]  $.\\
\indent  Division: $ A/B =[a^{l}, a^{u}]/[b^{l}, b^{u}]=[a^{l}, a^{u}][ \frac{1}{b^{l}}, \frac{1}{b^{u}}] $.
\end{Definition}

\begin{Definition}
(Interval-valued function) Let $ a > 0 $, $ b > 0 $. If the interval is of the from $ [a, b] $, the interval-valued function is take as $ h(k) = a^{(1-k)}b^{k} $ for $ k \in [0, 1] $.	
\end{Definition}
Let $ \hat{D}, \hat{m_{1}}, \hat{\delta_{1}}, \hat{\sigma_{1}}, \hat{m_{2}}, \hat{\delta_{2}}, \hat{\sigma_{2}}, \hat{\sigma_{3}}$ represent the interval numbers of $ D, {m_{1}}, {\delta_{1}}, {\sigma_{1}}, {m_{2}}, {\delta_{2}}, {\sigma_{2}}, {\sigma_{3}}$, respectively. The system (\ref{2.1}) with imprecise parameters becomes:
\begin{equation}\label{2.2}
\left\{
\begin{array}{lcl}
dS(t) = [\hat{D}(S^{0}-S(t))-\dfrac{\hat{m_{1}}S(t)x(t)}{\hat{\delta_{1}}}]dt + \hat{\sigma_{1}}S(t)dB_{1}(t) +  \int_{\mathbb{Y}}\gamma_{1}S(t^{-})\widetilde{N}(dt,du),\\
dx(t) = [\hat{m_{1}}S(t)x(t)-\hat{D}x(t)-\dfrac{\hat{m_{2}}x(t)y(t)}{\hat{\delta_{2}}}]dt + \hat{\sigma_{2}}S(t)dB_{2}(t) +  \int_{\mathbb{Y}}\gamma_{2}x(t^{-})\widetilde{N}(dt,du),  \\
dy(t) = [\hat{m_{2}}x(t)y(t)-\hat{D}y(t)]dt + \hat{\sigma_{3}}y(t)dB_{3}(t) +  \int_{\mathbb{Y}}\gamma_{3}y(t^{-})\widetilde{N}(dt,du).
\end{array}  
\right.
\end{equation}

where $ \hat{D}= [D^{l}, D^{u}] $, $ \hat{m_{1}}= [m_{1}^{l}, m_{1}^{u}] $, $ \hat{\delta_{1}}= [\delta_{1}^{l}, \delta_{1}^{u}] $, $ \hat{\sigma_{1}}= [\sigma_{1}^{l}, \sigma_{1}^{u}] $, $ \hat{m_{2}}= [m_{2}^{l}, m_{2}^{u}] $, $ \hat{\delta_{2}}= [\delta_{2}^{l}, \delta_{2}^{u}] $, $ \hat{\sigma_{2}}= [\sigma_{2}^{l}, \sigma_{2}^{u}] $, $ \hat{\sigma_{3}}= [\sigma_{3}^{l}, \sigma_{3}^{u}] $.  
According to the Theorem 1 in Pal et al. \cite{20} and considering the interval-valued function $ f (p) = (f^{l})^{1-p}(f^{u})^{p} $ for interval $ \hat{f}= [f^{l}, f^{u}] $  for $ p \in [0, 1] $, we can prove that system (\ref{2.2}) is equivalent to the following system:
\begin{equation}\label{2.3}
\left\{
\begin{array}{lcl}
dS(t) = [(D^{l})^{1-p} (D^{u})^{p}(S^{0}-S(t))-\dfrac{(m_{1}^{l})^{1-p} (m_{1}^{u})^{p}S(t)x(t)}{(\delta_{1}^{l})^{1-p} (\delta_{1}^{u})^{p}}]dt + (\sigma_{1}^{l})^{1-p} (\sigma_{1}^{u})^{p}S(t)dB_{1}(t) \\
 \ \ \ \ \ \ \ \ \ \ \ \  +  \int_{\mathbb{Y}}\gamma_{1}S(t^{-})\widetilde{N}(dt,du),\\
dx(t) = [(m_{1}^{l})^{1-p} (m_{1}^{u})^{p}S(t)x(t)-(D^{l})^{1-p} (D^{u})^{p}x(t)-\dfrac{(m_{2}^{l})^{1-p} (m_{2}^{u})^{p}x(t)y(t)}{(\delta_{2}^{l})^{1-p} (\delta_{2}^{u})^{p}}]dt\\
  \ \ \ \ \ \ \ \ \ \ \ \ +  (\sigma_{2}^{l})^{1-p} (\sigma_{2}^{u})^{p}S(t)dB_{2}(t) +  \int_{\mathbb{Y}}\gamma_{2}x(t^{-})\widetilde{N}(dt,du),  \\
dy(t) = [(m_{2}^{l})^{1-p} (m_{2}^{u})^{p}x(t)y(t)-(D^{l})^{1-p} (D^{u})^{p}y(t)]dt +  (\sigma_{3}^{l})^{1-p} (\sigma_{3}^{u})^{p}y(t)dB_{3}(t) +  \int_{\mathbb{Y}}\gamma_{3}y(t^{-})\widetilde{N}(dt,du),
\end{array}  
\right.
\end{equation}
for $ p \in [0, 1] $.
The corresponding deterministic system to (\ref{2.3}) is
\begin{equation}\label{2.4}
\left\{
\begin{array}{lcl}
dS(t) = [(D^{l})^{1-p} (D^{u})^{p}(S^{0}-S(t))-\dfrac{(m_{1}^{l})^{1-p} (m_{1}^{u})^{p}S(t)x(t)}{(\delta_{1}^{l})^{1-p} (\delta_{1}^{u})^{p}}]dt,\\
dx(t) = [(m_{1}^{l})^{1-p} (m_{1}^{u})^{p}S(t)x(t)-(D^{l})^{1-p} (D^{u})^{p}x(t)-\dfrac{(m_{2}^{l})^{1-p} (m_{2}^{u})^{p}x(t)y(t)}{(\delta_{2}^{l})^{1-p} (\delta_{2}^{u})^{p}}]dt,  \\
dy(t) = [(m_{2}^{l})^{1-p} (m_{2}^{u})^{p}x(t)y(t)-(D^{l})^{1-p} (D^{u})^{p}y(t)]dt.
\end{array}  
\right.
\end{equation}
Assume that there is a positive constant $ c $ such that
\begin{equation}\label{2.5}
\int_{\mathbb{Y}} [ \ln (1+ \gamma_{i}(u) )]^{2}\lambda (du) \leq c, \  i=1,2,3.
\end{equation}
For convenience in the following investigation, let $ (\Omega, \mathcal{F}, \{\mathcal{F}_{t}\}_{t\geq 0}, \mathbb{P}) $ be a complete probability space with a filtration $ \{\mathcal{F}_{t}\}_{t\geq 0} $ satisfying the usual conditions and define
\[
{\beta}_{i}= \dfrac{1}{2} ({ (\sigma_{i}^{l})^{1-p} (\sigma_{i}^{u})^{p}})^{2} + \int_{\mathbb{Y}} [ \gamma_{i}(u)- \ln (1+ \gamma_{i}(u) )]\lambda (du) , \  i=1,2,3.
\]
and
\[
R_{0}^{s}=\dfrac{S^{0}{(m_{1}^{l})^{1-p} (m_{1}^{u})^{p}}}{{(D^{l})^{1-p} (D^{u})^{p}}+\beta_{2}}, 
\]
\[
R_{1}^{s}= \dfrac{S^{0}{(m_{1}^{l})^{1-p} (m_{1}^{u})^{p}}{(m_{2}^{l})^{1-p} (m_{2}^{u})^{p}}{(\delta_{1}^{l})^{1-p} (\delta_{1}^{u})^{p}}}{{(m_{2}^{l})^{1-p} (m_{2}^{u})^{p}}{(\delta_{1}^{l})^{1-p} (\delta_{1}^{u})^{p}}({(D^{l})^{1-p} (D^{u})^{p}}+\beta_{2}) + {(m_{1}^{l})^{1-p} (m_{1}^{u})^{p}}({(D^{l})^{1-p} (D^{u})^{p}}+\beta_{3})}.
\]

\section{ Existence and uniqueness of positive solution of system (\ref{2.3}) }
\label{sec:3}
Assume for each $ m > 0 $ there exists $ L_{m} > 0 $ such that
\begin{description}
	\item[(H1)] $ \int_{\mathbb{Y}} |H_{i}(z,u)- H_{i}(\omega, u)|^{2}\lambda (du)\leq L_{m}|z-\omega|^{2} $, $ i=1,2,3, $ where $ H_{1}(z,u)= \gamma_{1}(u)z(t^{-})$, $ H_{2}(z,u)= \gamma_{2}(u)z(t^{-})$, $ H_{3}(z,u)= \gamma_{3}(u)z(t^{-})$ with $ |z|\vee |\omega|\leq m $.
	\item[(H2)] $ |\ln (1+\gamma_{i}(u))|\leq K_{i} $, for $ \gamma_{i}(u)>-1 $, where $ K_{i} $, $ i=1,2,3 $ are positive constants.
\end{description}

\begin{Theorem}
	Let Assumptions (H1) and (H2) hold. Then for any given initial value $ (S(0), x(0), y(0)) \in \mathbb{R}^{3}_{+} $, system (\ref{2.3}) has a unique solution $ (S(t), x(t), y(t)) \in \mathbb{R}^{3}_{+} $ for all $ t \geq 0 $ almost surely (a.s.).
\end{Theorem}

\noindent \textbf{Proof.} Because the coefficients of system (\ref{2.3}) are local Lipschitz continuous (Mao \cite{25}), for any given initial value $ (S(0), x(0), y(0)) \in \mathbb{R}^{3}_{+} $, there is a unique local solution $ (S(t), x(t), y(t)) $ on $ t\in [0,\tau_{e}) $, where $ \tau_{e} $ is the explosion time (see Mao \cite{25}). In order to show that the local solution is global, we only need to prove that $ \tau_{e}=\infty $ a.s. In this context, choosing a sufficiently large number $ k_{0}\geq 1 $ such that $ (S(0), x(0), y(0)) $ lie within the interval $ [\dfrac{1}{k_{0}}, k_{0}] $. For each integer $ k \geq k_{0} $, we define the stopping time as
\[
\tau_{k}=\inf \Big\{t\in [0, \tau_{e}): \min\{ S(t), x(t), y(t)\}\leq \dfrac{1}{k}   \textrm{ or }   \max \{ S(t), x(t), y(t)\}\geq k        \Big\}
\]
where $ \inf \emptyset= \infty $ ($\emptyset  $ being empty set). By the definition, $ \tau_{k} $ increases as $ k \rightarrow \infty $. Set $ \tau_{\infty}= \lim_{k \rightarrow \infty} \tau_{k} $. Hence $\tau_{\infty}\leq \tau_{e}  $ a.s. If $ \tau_{\infty}= \infty $ a.s. is true, then $ \tau_{e}=\infty  $ a.s. for all $ t>0 $. In other words, we need to verify $ \tau_{\infty}= \infty $ a.s. If this claim is wrong, then there exist a constant $ T > 0 $ and an $ \epsilon \in (0, 1) $ such that 
\[
\mathbb{P}\{ \tau_{\infty}\leq T \} > \epsilon.
\]
Hence there is an interger $ k_{1}\geq k_{0} $ such that
\begin{equation}\label{3.1}
\mathbb{P}\{ \tau_{k}\leq T \} \geq \epsilon \textrm{ for all }  k\geq k_{1} .
\end{equation}
Consider the Lyapunov function $ V:\mathbb{R}^{3}_{+} \rightarrow \overline{\mathbb{R}}_{+} $ defined for $ (S(t), x(t), y(t)) \in \mathbb{R}^{3}_{+} $  by
\begin{equation*}
\begin{split}
V(S,x,y)=& (S-A-A\ln\dfrac{S}{A}) + \dfrac{1}{{(\delta_{1}^{l})^{1-p} (\delta_{1}^{u})^{p}}}(x-B-B\ln \dfrac{x}{B}) +\\ &\dfrac{1}{{(\delta_{1}^{l})^{1-p} (\delta_{1}^{u})^{p}}{(\delta_{2}^{l})^{1-p} (\delta_{2}^{u})^{p}}}(y-1-\ln y),
\end{split}
\end{equation*}
where $ A,B $ are positive constants to be determined later. The nonnegativity of this function can be seen from
$ u-1-\ln u \geq 0 $ for any $ u>0 $.
Let $ k\geq k_{1} $ and $ T>0 $. Then, for any $ 0\leq t \leq \min{\tau_{k}, T} $, the It\v{o}’s formula (Situ \cite{26}) shows that
\begin{eqnarray*}
	dV (S, x, y) =& & LV(S, x, y)dt + (S- A){(\sigma_{1}^{l})^{1-p} (\sigma_{1}^{u})^{p}} dB_{1}(t)\\
	&& + \dfrac{1}{{(\delta_{1}^{l})^{1-p} (\delta_{1}^{u})^{p}}}(x-B){(\sigma_{2}^{l})^{1-p} (\sigma_{2}^{u})^{p}} dB_{2}(t)\\
	&& + \dfrac{1}{{(\delta_{1}^{l})^{1-p} (\delta_{1}^{u})^{p}}{(\delta_{2}^{l})^{1-p} (\delta_{2}^{u})^{p}}}(y-1){(\sigma_{3}^{l})^{1-p} (\sigma_{3}^{u})^{p}} dB_{3}(t)\\
	&& + \int_{\mathbb{Y}}[\gamma_{1}(u)S - A\ln (1+\gamma_{1}(u))]\widetilde{N}(dt, du)\\
	&& + \int_{\mathbb{Y}}[ \dfrac{1}{{(\delta_{1}^{l})^{1-p} (\delta_{1}^{u})^{p}}}(\gamma_{2}(u)x- B\ln (1+{(\delta_{2}^{l})^{1-p} (\delta_{2}^{u})^{p}}(u))  )   ]\widetilde{N}(dt, du)\\
	&& + \int_{\mathbb{Y}}[ \dfrac{1}{{(\delta_{1}^{l})^{1-p} (\delta_{1}^{u})^{p}}{(\delta_{2}^{l})^{1-p} (\delta_{2}^{u})^{p}}}( \gamma_{3}(u)y- \ln (1+ \gamma_{3}(u))  )  ]\widetilde{N}(dt, du),
\end{eqnarray*}
where $ L $ is a differential operator, and
\begin{eqnarray*}
	LV(S, x, y) &=& (1-\dfrac{A}{S}) [ {(D^{l})^{1-p} (D^{u})^{p}}  (S^{0}-S) - \dfrac{{(m_{1}^{l})^{1-p} (m_{1}^{u})^{p}}Sx}{{(\delta_{1}^{l})^{1-p} (\delta_{1}^{u})^{p}}}]\\
	&&  + \dfrac{A{(\sigma_{1}^{l})^{1-p} (\sigma_{1}^{u})^{p}}^{2}}{2}+ \int_{\mathbb{Y}}[A\gamma_{1}(u) - A\ln (1+\gamma_{1}(u))]\lambda (du)\\
	&& +  \dfrac{1}{{(\delta_{1}^{l})^{1-p} (\delta_{1}^{u})^{p}}} (1- \dfrac{B}{x}) ( {(m_{1}^{l})^{1-p} (m_{1}^{u})^{p}}Sx-{(D^{l})^{1-p} (D^{u})^{p}}x\\
	&&	-\dfrac{{(m_{2}^{l})^{1-p} (m_{2}^{u})^{p}}xy}{{(\delta_{2}^{l})^{1-p} (\delta_{2}^{u})^{p}}}) +\dfrac{B{(\sigma_{2}^{l})^{1-p} (\sigma_{2}^{u})^{p}}^{2}}{2{(\delta_{1}^{l})^{1-p} (\delta_{1}^{u})^{p}}}\\
	&& + \int_{\mathbb{Y}}[ \dfrac{1}{{(\delta_{1}^{l})^{1-p} (\delta_{1}^{u})^{p}}}(B\gamma_{2}(u)- B\ln (1+\gamma_{2}(u))) ]\lambda (du)\\
	&& +  \dfrac{1}{{(\delta_{1}^{l})^{1-p} (\delta_{1}^{u})^{p}}{(\delta_{2}^{l})^{1-p} (\delta_{2}^{u})^{p}}} (1- \dfrac{1}{y}) ( {(m_{2}^{l})^{1-p} (m_{2}^{u})^{p}}xy\\
	&&-{(D^{l})^{1-p} (D^{u})^{p}}y) + \dfrac{{(\sigma_{3}^{l})^{1-p} (\sigma_{3}^{u})^{p}}^{2}}{2{(\delta_{1}^{l})^{1-p} (\delta_{1}^{u})^{p}}{(\delta_{2}^{l})^{1-p} (\delta_{2}^{u})^{p}}}\\
	&& + \int_{\mathbb{Y}}[ \dfrac{1}{{(\delta_{1}^{l})^{1-p} (\delta_{1}^{u})^{p}}{(\delta_{2}^{l})^{1-p} (\delta_{2}^{u})^{p}}}( \gamma_{3}(u)- \ln (1+ \gamma_{3}(u))  )  ]\lambda (du)\\
	&=& {(D^{l})^{1-p} (D^{u})^{p}}(S^{0}-S) - \dfrac{A{(D^{l})^{1-p} (D^{u})^{p}}}{S}(S^{0}-S) \\
	&&+ (\dfrac{A{(m_{1}^{l})^{1-p} (m_{1}^{u})^{p}}}{{(\delta_{1}^{l})^{1-p} (\delta_{1}^{u})^{p}}}- \dfrac{{(D^{l})^{1-p} (D^{u})^{p}}}{{(\delta_{1}^{l})^{1-p} (\delta_{1}^{u})^{p}}} )x - \dfrac{B{(m_{1}^{l})^{1-p} (m_{1}^{u})^{p}}S}{{(\delta_{1}^{l})^{1-p} (\delta_{1}^{u})^{p}}}\\
	&& + \dfrac{B{(D^{l})^{1-p} (D^{u})^{p}}}{{(\delta_{1}^{l})^{1-p} (\delta_{1}^{u})^{p}}} +  (\dfrac{B{(m_{2}^{l})^{1-p} (m_{2}^{u})^{p}}}{{(\delta_{1}^{l})^{1-p} (\delta_{1}^{u})^{p}}{(\delta_{2}^{l})^{1-p} (\delta_{2}^{u})^{p}}}\\
	&&- \dfrac{{(D^{l})^{1-p} (D^{u})^{p}}}{{(\delta_{1}^{l})^{1-p} (\delta_{1}^{u})^{p}}{(\delta_{2}^{l})^{1-p} (\delta_{2}^{u})^{p}}} )y - \dfrac{{(m_{2}^{l})^{1-p} (m_{2}^{u})^{p}}x}{{(\delta_{1}^{l})^{1-p} (\delta_{1}^{u})^{p}}{(\delta_{2}^{l})^{1-p} (\delta_{2}^{u})^{p}}} + \dfrac{{(D^{l})^{1-p} (D^{u})^{p}}}{{(\delta_{1}^{l})^{1-p} (\delta_{1}^{u})^{p}}{(\delta_{2}^{l})^{1-p} (\delta_{2}^{u})^{p}}}\\
	&& + \dfrac{A{(\sigma_{1}^{l})^{1-p} (\sigma_{1}^{u})^{p}}^{2}}{2} + \dfrac{B{(\sigma_{2}^{l})^{1-p} (\sigma_{2}^{u})^{p}}^{2}}{2{(\delta_{1}^{l})^{1-p} (\delta_{1}^{u})^{p}}} + \dfrac{{(\sigma_{3}^{l})^{1-p} (\sigma_{3}^{u})^{p}}^{2}}{2{(\delta_{1}^{l})^{1-p} (\delta_{1}^{u})^{p}}{(\delta_{2}^{l})^{1-p} (\delta_{2}^{u})^{p}}}\\
	&& + A\int_{\mathbb{Y}}[\gamma_{1}(u) - \ln (1+\gamma_{1}(u))]\lambda (du)
\end{eqnarray*}
\begin{eqnarray*}
	&& + \dfrac{B}{{(\delta_{1}^{l})^{1-p} (\delta_{1}^{u})^{p}}}\int_{\mathbb{Y}}[ \gamma_{2}(u)- \ln (1+\gamma_{2}(u)) ]\lambda (du)\\
	&& + \dfrac{1}{{(\delta_{1}^{l})^{1-p} (\delta_{1}^{u})^{p}}{(\delta_{2}^{l})^{1-p} (\delta_{2}^{u})^{p}}} \int_{\mathbb{Y}}[ \gamma_{3}(u)- \ln (1+ \gamma_{3}(u))   ]\lambda (du) .
\end{eqnarray*}
Choose $ A= \dfrac{{(D^{l})^{1-p} (D^{u})^{p}}}{{(m_{1}^{l})^{1-p} (m_{1}^{u})^{p}}} $, $ B= \dfrac{{(D^{l})^{1-p} (D^{u})^{p}}}{{(m_{2}^{l})^{1-p} (m_{2}^{u})^{p}}} $, by use of inequation $x-\ln(x+1)\geq 0 $ for $ x>-1 $ and (H2), 
\begin{eqnarray*}
	LV(S, x, y) &=& {(D^{l})^{1-p} (D^{u})^{p}}(S^{0}-S) - \dfrac{A{(D^{l})^{1-p} (D^{u})^{p}}}{S}(S^{0}-S) - \dfrac{B{(m_{1}^{l})^{1-p} (m_{1}^{u})^{p}}S}{{(\delta_{1}^{l})^{1-p} (\delta_{1}^{u})^{p}}}\\
	&& + \dfrac{B{(D^{l})^{1-p} (D^{u})^{p}}}{{(\delta_{1}^{l})^{1-p} (\delta_{1}^{u})^{p}}}  - \dfrac{{(m_{2}^{l})^{1-p} (m_{2}^{u})^{p}}x}{{(\delta_{1}^{l})^{1-p} (\delta_{1}^{u})^{p}}{(\delta_{2}^{l})^{1-p} (\delta_{2}^{u})^{p}}} + \dfrac{{(D^{l})^{1-p} (D^{u})^{p}}}{{(\delta_{1}^{l})^{1-p} (\delta_{1}^{u})^{p}}{(\delta_{2}^{l})^{1-p} (\delta_{2}^{u})^{p}}}\\
	&&+ \dfrac{A({(\sigma_{1}^{l})^{1-p} (\sigma_{1}^{u})^{p}})^{2}}{2} + \dfrac{B({(\sigma_{2}^{l})^{1-p} (\sigma_{2}^{u})^{p}})^{2}}{2{(\delta_{1}^{l})^{1-p} (\delta_{1}^{u})^{p}}} + \dfrac{({(\sigma_{3}^{l})^{1-p} (\sigma_{3}^{u})^{p}})^{2}}{2{(\delta_{1}^{l})^{1-p} (\delta_{1}^{u})^{p}}{(\delta_{2}^{l})^{1-p} (\delta_{2}^{u})^{p}}}\\
	&& + A\int_{\mathbb{Y}}[\gamma_{1}(u) - \ln (1+\gamma_{1}(u))]\lambda (du)\\
	&& + \dfrac{B}{{(\delta_{1}^{l})^{1-p} (\delta_{1}^{u})^{p}}}\int_{\mathbb{Y}}[ \gamma_{2}(u)- \ln (1+\gamma_{2}(u)) ]\lambda (du)\\
	&& + \dfrac{1}{{(\delta_{1}^{l})^{1-p} (\delta_{1}^{u})^{p}}{(\delta_{2}^{l})^{1-p} (\delta_{2}^{u})^{p}}} \int_{\mathbb{Y}}[ \gamma_{3}(u)- \ln (1+ \gamma_{3}(u))   ]\lambda (du) \\
	& \leq & {(D^{l})^{1-p} (D^{u})^{p}} S^{0} + \dfrac{({(D^{l})^{1-p} (D^{u})^{p}})^{2}}{{(m_{1}^{l})^{1-p} (m_{1}^{u})^{p}}} + \dfrac{({(D^{l})^{1-p} (D^{u})^{p}})^{2}}{{(m_{2}^{l})^{1-p} (m_{2}^{u})^{p}}{(\delta_{1}^{l})^{1-p} (\delta_{1}^{u})^{p}}}
	\\
	&&+ \dfrac{{(D^{l})^{1-p} (D^{u})^{p}}}{{(\delta_{1}^{l})^{1-p} (\delta_{1}^{u})^{p}}{(\delta_{2}^{l})^{1-p} (\delta_{2}^{u})^{p}}}
	 + \dfrac{{(D^{l})^{1-p} (D^{u})^{p}}({(\sigma_{1}^{l})^{1-p} (\sigma_{1}^{u})^{p}})^{2}}{2{(m_{1}^{l})^{1-p} (m_{1}^{u})^{p}}} \\
	 &&+ \dfrac{{(D^{l})^{1-p} (D^{u})^{p}}({(\sigma_{2}^{l})^{1-p} (\sigma_{2}^{u})^{p}})^{2}}{2{(m_{2}^{l})^{1-p} (m_{2}^{u})^{p}}{(\delta_{1}^{l})^{1-p} (\delta_{1}^{u})^{p}}} + \dfrac{({(\sigma_{3}^{l})^{1-p} (\sigma_{3}^{u})^{p}})^{2}}{2{(\delta_{1}^{l})^{1-p} (\delta_{1}^{u})^{p}}{(\delta_{2}^{l})^{1-p} (\delta_{2}^{u})^{p}}} + 3C \\
	&:=& \widetilde{K},
\end{eqnarray*}
where $\widetilde{K}  $ is a positive constant and
 \begin{eqnarray*} C&=& \max \big\{   A\int_{\mathbb{Y}}[\gamma_{1}(u) - \ln (1+\gamma_{1}(u))]\lambda (du), \frac{B}{{(\delta_{1}^{l})^{1-p} (\delta_{1}^{u})^{p}}}\int_{\mathbb{Y}}[ \gamma_{2}(u)\\
 	&&- \ln (1+\gamma_{2}(u)) ]\lambda (du), \frac{1}{{(\delta_{1}^{l})^{1-p} (\delta_{1}^{u})^{p}}{(\delta_{2}^{l})^{1-p} (\delta_{2}^{u})^{p}}} \int_{\mathbb{Y}}[ \gamma_{3}(u)- \ln (1+ \gamma_{3}(u))   ]\lambda (du)           \big\} .\end{eqnarray*}
Therefore,
\begin{equation}\label{3.2}
\begin{split}
\int_{0}^{\tau_{k} \wedge T} dV(S, x, y) \leq& \int_{0}^{\tau_{k} \wedge T}\widetilde{K}ds + \int_{0}^{\tau_{k} \wedge T} (S-A){(\sigma_{1}^{l})^{1-p} (\sigma_{1}^{u})^{p}}dB_{1}(s) \\
& +\int_{0}^{\tau_{k} \wedge T} \dfrac{1}{{(\delta_{1}^{l})^{1-p} (\delta_{1}^{u})^{p}}}(x-B){(\sigma_{2}^{l})^{1-p} (\sigma_{2}^{u})^{p}}dB_{2}(s) \\
&+ \dfrac{1}{{(\delta_{1}^{l})^{1-p} (\delta_{1}^{u})^{p}}{(\delta_{2}^{l})^{1-p} (\delta_{2}^{u})^{p}}}(y-1){(\sigma_{3}^{l})^{1-p} (\sigma_{3}^{u})^{p}}dB_{3}(s) \\
& + \int_{0}^{\tau_{k} \wedge T} \int_{\mathbb{Y}}[\gamma_{1}(u)S - A\ln (1+\gamma_{1}(u))]\widetilde{N}(ds, du)\\
& +  \int_{0}^{\tau_{k} \wedge T}\int_{\mathbb{Y}}[ \dfrac{1}{{(\delta_{1}^{l})^{1-p} (\delta_{1}^{u})^{p}}}(\gamma_{2}(u)x- B\ln (1+{(\delta_{2}^{l})^{1-p} (\delta_{2}^{u})^{p}}(u))  )   ]\widetilde{N}(ds, du)\\
& +  \int_{0}^{\tau_{k} \wedge T} \int_{\mathbb{Y}}[ \dfrac{1}{{(\delta_{1}^{l})^{1-p} (\delta_{1}^{u})^{p}}{(\delta_{2}^{l})^{1-p} (\delta_{2}^{u})^{p}}}( \gamma_{3}(u)y- \ln (1+ \gamma_{3}(u))  )  ]\widetilde{N}(ds, du).
\end{split}
\end{equation}

Taking the expectations on both sides of (\ref{3.2}), we obtain that
\begin{equation}\label{3.3}
\begin{split}
\mathbb{E}V(S(\tau_{k} \wedge T ), x(\tau_{k} \wedge T ), y(\tau_{k} \wedge T) )
&\leq V( S(0), x(0), y(0) ) + \widetilde{K} \mathbb{E}(\tau_{k} \wedge T )\\
&\leq V( S(0), x(0), y(0) ) + \widetilde{K}T.
\end{split}
\end{equation}
Let $ \Omega_{k}=\{ \tau_{k} \wedge T \} $ for $ k\geq k_{1} $. Then, by (\ref{3.1}), we know that $ \mathbb{P}(\Omega_{k})\geq \epsilon$. Noting that for every $ \omega\in\Omega_{k}  $, there exist $ S(\tau_{k}, \omega ), x(\tau_{k}, \omega)$ and $ y(\tau_{k}, \omega) $, all of which equal either $ k $ or $ \frac{1}{k} $. Hence $V(S(\tau_{k}, \omega ), x(\tau_{k}, \omega), y(\tau_{k}, \omega))  $ is no less than
$ k-A-A\ln\frac{k}{A} $ or $ \frac{1}{k}-A-A\ln \frac{1}{Ak} $ or $ \frac{1}{{(\delta_{1}^{l})^{1-p} (\delta_{1}^{u})^{p}}}(k-B-B\ln\frac{k}{B}) $ or $ \frac{1}{{(\delta_{1}^{l})^{1-p} (\delta_{1}^{u})^{p}}}(\frac{1}{k}-B-B\ln\frac{1}{Bk}) $ or $ \frac{1}{{(\delta_{1}^{l})^{1-p} (\delta_{1}^{u})^{p}}{(\delta_{2}^{l})^{1-p} (\delta_{2}^{u})^{p}}}(k-1-\ln k) $ or $ \frac{1}{{(\delta_{1}^{l})^{1-p} (\delta_{1}^{u})^{p}}{(\delta_{2}^{l})^{1-p} (\delta_{2}^{u})^{p}}}(\frac{1}{k}-1+ \ln k) $.
Hence
\begin{equation*}
\begin{split}
M:= & (k-A-A\ln\frac{k}{A}) \wedge (\frac{1}{k}-A-A\ln \frac{1}{Ak}) \wedge (\frac{1}{{(\delta_{1}^{l})^{1-p} (\delta_{1}^{u})^{p}}}(k-B-B\ln\frac{k}{B}))\\  
&\wedge (\frac{1}{{(\delta_{1}^{l})^{1-p} (\delta_{1}^{u})^{p}}}(\frac{1}{k}-B-B\ln\frac{1}{Bk})) \wedge (\frac{1}{{(\delta_{1}^{l})^{1-p} (\delta_{1}^{u})^{p}}{(\delta_{2}^{l})^{1-p} (\delta_{2}^{u})^{p}}}(k-1-\ln k))\\
& \wedge (\frac{1}{{(\delta_{1}^{l})^{1-p} (\delta_{1}^{u})^{p}}{(\delta_{2}^{l})^{1-p} (\delta_{2}^{u})^{p}}}(\frac{1}{k}-1+ \ln k) ).
\end{split}
\end{equation*}
Thus, by (\ref{3.3}), we konw that
\[
V(S(0), x(0), y(0) ) + \widetilde{K}T \geq \mathbb{E}[1_{\Omega_{k}(\omega)} V(S(\tau_{k}, \omega ), x(\tau_{k}, \omega), y(\tau_{k}, \omega)) ] \geq \epsilon M,
\]
where $ 1_{\Omega_{k}(\omega)} $ represents the indicator function of $ \Omega_{k}(\omega) $. Setting $ k \rightarrow \infty $ leads to the contradiction
\[
\infty > V ( S(0), x(0), y(0)  ) +  \widetilde{K}T = \infty.
\]
Therefore, we have  $ \tau_{\infty}= \infty $ a.s. The proof is complete. \qed

\section{Extinction}
\label{sec:4}

When studying dynamical models, two of the most interesting issues are persistence and extinction. In this
section, we discuss the extinction of microorganisms in system (\ref{2.3}) and leave its persistence to the next section.

Before we state the main results of this section, several lemmas (gao et al. \cite{27}) should be recalled without the proofs and relevant explanations.

\begin{description}
	\item[(H3)] 	Assume that for some $ \theta >2 $, $ {(D^{l})^{1-p} (D^{u})^{p}} - \dfrac{\theta-1}{2}\sigma^{2} - \dfrac{\zeta}{\theta} > 0$, where $\zeta= \int_{\mathbb{Y}}[(1+(\gamma_{1}(u)\vee \gamma_{2}(u) \vee \gamma_{3}(u)))^{p} -1 - \gamma_{1}(u)\wedge \gamma_{2}(u) \wedge \gamma_{3}(u)  ]\lambda (du)  $, $ \sigma^{2} = ({(\sigma_{1}^{l})^{1-p} (\sigma_{1}^{u})^{p}})^{2} \vee ({(\sigma_{2}^{l})^{1-p} (\sigma_{2}^{u})^{p}})^{2} \vee ({(\sigma_{3}^{l})^{1-p} (\sigma_{3}^{u})^{p}})^{2} $.
\end{description}

\begin{Lemma}\label{L41}
	Let Assumption (H3) hold.  For any initial value $ (S(0), x(0), y(0)) \in \mathbb{R}^{3}_{+} $, there is a unique positive solution $ (S(t), x(t), y(t)) $ to system (\ref{2.3}) and the solution will remain in $\mathbb{R}^{3}_{+} $ with probability one, i.e., $ (S(t), x(t), y(t))\in \mathbb{R}^{3}_{+} $ for all $ t \geq 0 $ a.s. Then
	\[
	\lim_{t\rightarrow \infty}\dfrac{S(t)+ x(t)+ y(t)}{t}=0 \ \ a.s.
	\]
	
	\[
	\lim_{t\rightarrow \infty}\dfrac{S(t)}{t}=0, \ \ 	\lim_{t\rightarrow \infty}\dfrac{x(t)}{t}=0, \ \ 	\lim_{t\rightarrow \infty}\dfrac{y(t)}{t}=0 \ \ a.s.
	\]
	Moreover
	
	\[
	\lim_{t\rightarrow \infty}\dfrac{\int_{0}^{t}S(s)dB_{1}(s)}{t} =0, \ 	\lim_{t\rightarrow \infty}\dfrac{\int_{0}^{t}x(s)dB_{2}(s)}{t} =0, \ 	\lim_{t\rightarrow \infty}\dfrac{\int_{0}^{t}y(s)dB_{3}(s)}{t} =0 \ \ a.s.
	\]
	\[
	\lim_{t\rightarrow \infty}\dfrac{\int_{0}^{t}\int_{\mathbb{Y}}\gamma_{1}(u)S(s^{-})\widetilde{N}(ds,du)}{t} =0, \  	\lim_{t\rightarrow \infty}\dfrac{\int_{0}^{t}\int_{\mathbb{Y}}\gamma_{2}(u)x(s^{-})\widetilde{N}(ds,du)}{t} =0,
	\]
	\[
	\lim_{t\rightarrow \infty}\dfrac{\int_{0}^{t}\int_{\mathbb{Y}}\gamma_{3}(u)y(s^{-})\widetilde{N}(ds,du)}{t} =0 \ a.s.
	\]
\end{Lemma}

\begin{Lemma}\label{L42}
	For any given initial value $ (S(0), x(0), y(0)) \in \mathbb{R}^{3}_{+} $， the solution $ (S(t), x(t), y(t)) $ of system (\ref{2.3})  has the following	property
	\[
	\limsup_{t\rightarrow \infty}\dfrac{\ln x(t)}{t}\leq 0, \ \ \limsup_{t\rightarrow \infty}\dfrac{\ln y(t)}{t}\leq 0 \ \ a.s.
	\]
\end{Lemma}

\begin{Theorem}
	Let Assumption (H3) hold.  For any initial value $ (S(0), x(0), y(0)) \in \mathbb{R}^{3}_{+} $, there is a unique positive solution $ (S(t), x(t), y(t)) $ to system (\ref{2.3}) and the solution will remain in $\mathbb{R}^{3}_{+} $ with probability one, i.e., $ (S(t), x(t), y(t))\in \mathbb{R}^{3}_{+} $ for all $ t \geq 0 $ a.s. If $ R_{0}^{s} <1$, then
	\[
	\limsup_{t\rightarrow \infty}\dfrac{\ln x(t)}{t} \leq ({(D^{l})^{1-p} (D^{u})^{p}} + {\beta}_{2})(R_{0}^{s} - 1) < 0 \ \ a.s.
	\]
	and
	\[
	\limsup_{t\rightarrow \infty}\dfrac{\ln y(t)}{t} \leq -({(D^{l})^{1-p} (D^{u})^{p}} + {\beta}_{3}) < 0 \ \ a.s.
	\]
	Moreover
	\[
	\lim_{t\rightarrow \infty} \langle S \rangle_{t} = S^{0} \ \  a.s.
	\]
	That is to say, microbes $ x $ and $ y $ will become extinct exponentially with probability one.
\end{Theorem}

\noindent \textbf{Proof.}
Integrating from $ 0 $ to $ t $ on both sides of (\ref{2.3}), yields
\begin{equation*}
\begin{split}
&\dfrac{S(t)-S(0)}{t} + \dfrac{x(t)-x(0)}{{(\delta_{1}^{l})^{1-p} (\delta_{1}^{u})^{p}}t} + \dfrac{y(t)-y(0)}{{(\delta_{1}^{l})^{1-p} (\delta_{1}^{u})^{p}}{(\delta_{2}^{l})^{1-p} (\delta_{2}^{u})^{p}}t}\\
= & {(D^{l})^{1-p} (D^{u})^{p}}S^{0} - {(D^{l})^{1-p} (D^{u})^{p}}\langle S \rangle_{t} - \dfrac{{(D^{l})^{1-p} (D^{u})^{p}}}{{(\delta_{1}^{l})^{1-p} (\delta_{1}^{u})^{p}}}\langle x \rangle_{t}\\
& - \dfrac{{(D^{l})^{1-p} (D^{u})^{p}}}{{(\delta_{1}^{l})^{1-p} (\delta_{1}^{u})^{p}}{(\delta_{2}^{l})^{1-p} (\delta_{2}^{u})^{p}}}\langle y \rangle_{t}  
+ \dfrac{{(\sigma_{1}^{l})^{1-p} (\sigma_{1}^{u})^{p}}}{t} \int_{0}^{t} S(s) dB_{1}(s)\\
& + \dfrac{1}{t} \int_{0}^{t}\int_{\mathbb{Y}}\gamma_{1}(u)S(s^{-})\widetilde{N}(ds,du)
+ \dfrac{{(\sigma_{2}^{l})^{1-p} (\sigma_{2}^{u})^{p}}}{{(\delta_{1}^{l})^{1-p} (\delta_{1}^{u})^{p}}t} \int_{0}^{t} x(s) dB_{2}(s) \\
&+ \dfrac{1}{{(\delta_{1}^{l})^{1-p} (\delta_{1}^{u})^{p}}t} \int_{0}^{t}\int_{\mathbb{Y}}\gamma_{2}(u)x(s^{-})\widetilde{N}(ds,du)\\
&+ \dfrac{{(\sigma_{3}^{l})^{1-p} (\sigma_{3}^{u})^{p}}}{{(\delta_{1}^{l})^{1-p} (\delta_{1}^{u})^{p}}{(\delta_{2}^{l})^{1-p} (\delta_{2}^{u})^{p}}t} \int_{0}^{t} y(s) dB_{3}(s)\\
& + \dfrac{1}{{(\delta_{1}^{l})^{1-p} (\delta_{1}^{u})^{p}}{(\delta_{2}^{l})^{1-p} (\delta_{2}^{u})^{p}}t} \int_{0}^{t}\int_{\mathbb{Y}}\gamma_{3}(u)y(s^{-})\widetilde{N}(ds,du).
\end{split}
\end{equation*}
Clearly, we can derive that
\begin{equation}\label{4.1}
\langle S \rangle_{t} = S^{0} - \dfrac{1}{{(\delta_{1}^{l})^{1-p} (\delta_{1}^{u})^{p}}} \langle x \rangle_{t} - \dfrac{1}{{(\delta_{1}^{l})^{1-p} (\delta_{1}^{u})^{p}}{(\delta_{2}^{l})^{1-p} (\delta_{2}^{u})^{p}}}\langle y \rangle_{t} + \phi (t),
\end{equation}
where 
\begin{equation*}
\begin{split}
\phi (t) = & \dfrac{1}{{(D^{l})^{1-p} (D^{u})^{p}}}[- \dfrac{S(t)-S(0)}{t} - \dfrac{x(t)-x(0)}{{(\delta_{1}^{l})^{1-p} (\delta_{1}^{u})^{p}}t} \\
&- \dfrac{y(t)-y(0)}{{(\delta_{1}^{l})^{1-p} (\delta_{1}^{u})^{p}}{(\delta_{2}^{l})^{1-p} (\delta_{2}^{u})^{p}}t}   
+ \dfrac{{(\sigma_{1}^{l})^{1-p} (\sigma_{1}^{u})^{p}}}{t} \int_{0}^{t} S(s) dB_{1}(s)\\
&+ \dfrac{1}{t} \int_{0}^{t}\int_{\mathbb{Y}}\gamma_{1}(u)S(s^{-})\widetilde{N}(ds,du)
+ \dfrac{{(\sigma_{2}^{l})^{1-p} (\sigma_{2}^{u})^{p}}}{{(\delta_{1}^{l})^{1-p} (\delta_{1}^{u})^{p}}t} \int_{0}^{t} x(s) dB_{2}(s)\\
&+ \dfrac{1}{{(\delta_{1}^{l})^{1-p} (\delta_{1}^{u})^{p}}t} \int_{0}^{t}\int_{\mathbb{Y}}\gamma_{2}(u)x(s^{-})\widetilde{N}(ds,du)\\
&+\dfrac{{(\sigma_{3}^{l})^{1-p} (\sigma_{3}^{u})^{p}}}{{(\delta_{1}^{l})^{1-p} (\delta_{1}^{u})^{p}}{(\delta_{2}^{l})^{1-p} (\delta_{2}^{u})^{p}}t} \int_{0}^{t} y(s) dB_{3}(s) \\
&+ \dfrac{1}{{(\delta_{1}^{l})^{1-p} (\delta_{1}^{u})^{p}}{(\delta_{2}^{l})^{1-p} (\delta_{2}^{u})^{p}}t} \int_{0}^{t}\int_{\mathbb{Y}}\gamma_{3}(u)y(s^{-})\widetilde{N}(ds,du)].
\end{split}
\end{equation*}
This together with Lemma~\ref{L41} implies
\begin{equation}\label{4.2}
\lim_{t\rightarrow \infty}\phi (t) = 0 \ \ a.s. 
\end{equation}
Applying Ito formula to (\ref{2.3}) we can conclude that
\begin{equation}\label{4.3}
\begin{split}
d\ln x =& [ {(m_{1}^{l})^{1-p} (m_{1}^{u})^{p}}S(t) - \dfrac{{(m_{2}^{l})^{1-p} (m_{2}^{u})^{p}}}{{(\delta_{2}^{l})^{1-p} (\delta_{2}^{u})^{p}}}y(t) - {(D^{l})^{1-p} (D^{u})^{p}} - \beta_{2} ] dt \\
&+ {(\sigma_{2}^{l})^{1-p} (\sigma_{2}^{u})^{p}} dB_{2}(t) + \int_{\mathbb{Y}}[\ln (1+\gamma_{2}(u) )] \widetilde{N}(dt,du).
\end{split}
\end{equation}
Integrating (\ref{4.3}) from $ 0 $ to $ t $ and then dividing by $ t $ on both sides, we obtain
\begin{equation}\label{4.4}
\begin{split}
\dfrac{\ln x(t)}{t} =& {(m_{1}^{l})^{1-p} (m_{1}^{u})^{p}} \langle S \rangle_{t} - \dfrac{{(m_{2}^{l})^{1-p} (m_{2}^{u})^{p}}}{{(\delta_{2}^{l})^{1-p} (\delta_{2}^{u})^{p}}}\langle y \rangle_{t} - {(D^{l})^{1-p} (D^{u})^{p}} - \beta_{2} \\
&+ \dfrac{{(m_{2}^{l})^{1-p} (m_{2}^{u})^{p}}(t)}{t} + \dfrac{\widetilde{M}_{2}(t)}{t} + \dfrac{\ln x(0)}{t}.
\end{split}
\end{equation}
where $ M_{i}(t): = \int_{0}^{t}{(\sigma_{i}^{l})^{1-p} (\sigma_{i}^{u})^{p}}dB_{i}(s) $ and $ \widetilde{M}_{i}(t):= \int_{0}^{t}\int_{\mathbb{Y}}[\ln (1+\gamma_{i}(u) )] \widetilde{N}(ds,du)   $, $ i=1,2,3 $ are all martingale terms.
Substituting (\ref{4.1}) into (\ref{4.4}) yields that
\begin{equation}\label{4.5}
\begin{split}
\dfrac{\ln x(t)}{t} =& {(m_{1}^{l})^{1-p} (m_{1}^{u})^{p}} [ S^{0} - \dfrac{1}{{(\delta_{1}^{l})^{1-p} (\delta_{1}^{u})^{p}}}\langle x \rangle_{t} -  \dfrac{1}{{(\delta_{1}^{l})^{1-p} (\delta_{1}^{u})^{p}}{(\delta_{2}^{l})^{1-p} (\delta_{2}^{u})^{p}}}\langle y \rangle_{t}\\
& + \phi(t)] - \dfrac{{(m_{2}^{l})^{1-p} (m_{2}^{u})^{p}}}{{(\delta_{2}^{l})^{1-p} (\delta_{2}^{u})^{p}}}\langle y \rangle_{t} - {(D^{l})^{1-p} (D^{u})^{p}} - \beta_{2}\\
& + \dfrac{{(m_{2}^{l})^{1-p} (m_{2}^{u})^{p}}(t)}{t} + \dfrac{\widetilde{M}_{2}(t)}{t} + \dfrac{\ln x(0)}{t}\\
=& {(m_{1}^{l})^{1-p} (m_{1}^{u})^{p}} S^{0} - {(D^{l})^{1-p} (D^{u})^{p}} - \beta_{2}- \dfrac{{(m_{1}^{l})^{1-p} (m_{1}^{u})^{p}}}{{(\delta_{1}^{l})^{1-p} (\delta_{1}^{u})^{p}}}\langle x \rangle_{t}\\
&-\dfrac{{(m_{1}^{l})^{1-p} (m_{1}^{u})^{p}}}{{(\delta_{1}^{l})^{1-p} (\delta_{1}^{u})^{p}}{(\delta_{2}^{l})^{1-p} (\delta_{2}^{u})^{p}}}\langle y \rangle_{t}-\dfrac{{(m_{2}^{l})^{1-p} (m_{2}^{u})^{p}}}{{(\delta_{2}^{l})^{1-p} (\delta_{2}^{u})^{p}}}\langle y \rangle_{t}\\
&+ {(m_{1}^{l})^{1-p} (m_{1}^{u})^{p}}\phi(t)+ \dfrac{{(m_{2}^{l})^{1-p} (m_{2}^{u})^{p}}(t)}{t} + \dfrac{\widetilde{M}_{2}(t)}{t} + \dfrac{\ln x(0)}{t}\\
\leq & {(m_{1}^{l})^{1-p} (m_{1}^{u})^{p}} S^{0} - {(D^{l})^{1-p} (D^{u})^{p}} - \beta_{2}+ {(m_{1}^{l})^{1-p} (m_{1}^{u})^{p}}\phi(t)\\
&+ \dfrac{{(m_{2}^{l})^{1-p} (m_{2}^{u})^{p}}(t)}{t} + \dfrac{\widetilde{M}_{2}(t)}{t} + \dfrac{\ln x(0)}{t}.
\end{split}
\end{equation}
Moreover, according to (\ref{2.5}),
\[
\langle M_{i}, M_{i} \rangle_{t} = \int_{0}^{t} ({(\sigma_{i}^{l})^{1-p} (\sigma_{i}^{u})^{p}})^{2} ds = ({(\sigma_{i}^{l})^{1-p} (\sigma_{i}^{u})^{p}})^{2} t,
\]
\[
\langle \widetilde{M}_{i}, \widetilde{M}_{i} \rangle_{t} = \int_{0}^{t} \int_{\mathbb{Y}}[\ln (1+\gamma_{i}(u) )]^{2} \lambda(du)ds \leq ct.
\]
where $ \langle M_{i}, M_{i} \rangle_{t} $ and $ \langle \widetilde{M}_{i}, \widetilde{M}_{i} \rangle_{t} $ represent the quadratic variation of $M_{i}$ and $ \widetilde{M}_{i} $, respectively. Thus, by strong law of large numbers, we have
\begin{equation}\label{4.6}
\lim_{t\rightarrow \infty} \dfrac{M_{i}(t)}{t} = 0 \ \ a.s. \textrm{ and } \lim_{t\rightarrow \infty} \dfrac{\widetilde{M}_{i}(t)}{t} = 0 \ \ a.s.
\end{equation}
Clearly, if $ R_{0}^{s} < 1$ holds, then taking the superior limit on both sides of (\ref{4.5}) and combining with (\ref{4.2}), (\ref{4.6}), we know that
\begin{equation*}
\begin{split}
\limsup_{t\rightarrow \infty}\dfrac{\ln x(t)}{t} \leq& {(m_{1}^{l})^{1-p} (m_{1}^{u})^{p}}S^{0} -{(D^{l})^{1-p} (D^{u})^{p}} - \beta_{2} \\
=& ( {(D^{l})^{1-p} (D^{u})^{p}} + \beta_{2})(R_{0}^{s}-1)\\
<&0 \ \ a.s.
\end{split}
\end{equation*}
which implies
\begin{equation}\label{4.7}
\lim_{t\rightarrow \infty} x(t) = 0 \ \ a.s.
\end{equation}
Similarly, we get
\begin{equation}\label{4.8}
\dfrac{\ln y(t)}{t} = {(m_{2}^{l})^{1-p} (m_{2}^{u})^{p}} \langle x \rangle_{t} - {(D^{l})^{1-p} (D^{u})^{p}} - \beta_{3} + \dfrac{M_{3}(t)}{t} + \dfrac{\widetilde{M}_{3}(t)}{t} + \dfrac{\ln y(0)}{t}.
\end{equation}
Thus,
\[
\limsup_{t\rightarrow \infty}\dfrac{\ln y(t)}{t}  = -( {(D^{l})^{1-p} (D^{u})^{p}} + \beta_{3})<0 \ \ a.s.
\]
which yields 
\begin{equation}\label{4.9}
\lim_{t\rightarrow \infty} y(t) = 0 \ \ a.s.
\end{equation}
It is now easy to derive from (\ref{4.8}) and (\ref{4.9}) that
\[
\lim_{t\rightarrow \infty} \langle S \rangle_{t} = S^{0} \ \  a.s.
\]
This completes the proof.  \qed

\begin{Theorem}
	Let Assumption (H3) hold. For any initial value $ (S(0), x(0), y(0))$ $ \in \mathbb{R}^{3}_{+} $, there is a unique positive solution $ (S(t), x(t), y(t)) $ to system (\ref{2.3}) and the solution will remain in $\mathbb{R}^{3}_{+} $ with probability one, i.e., $ (S(t), x(t), y(t))\in \mathbb{R}^{3}_{+} $ for all $ t \geq 0 $ a.s. If $ R_{1}^{s} <1<R_{0}^{s}$, then
	\begin{equation*}
	\begin{split}
	\limsup_{t\rightarrow \infty}\dfrac{\ln y(t)}{t} \leq& [\dfrac{{(m_{2}^{l})^{1-p} (m_{2}^{u})^{p}}{(\delta_{1}^{l})^{1-p} (\delta_{1}^{u})^{p}}}{{(m_{1}^{l})^{1-p} (m_{1}^{u})^{p}}}( {(D^{l})^{1-p} (D^{u})^{p}} + \beta_{2})  \\
	&+ {(D^{l})^{1-p} (D^{u})^{p}} + \beta_{3} ](R_{1}^{s} - 1) \\
	<& 0 \ \ a.s.
	\end{split}
	\end{equation*}
	Moreover	
	\[
	\lim_{t\rightarrow \infty} \langle S \rangle_{t} = \dfrac{ {(D^{l})^{1-p} (D^{u})^{p}} + \beta_{2}}{{(m_{1}^{l})^{1-p} (m_{1}^{u})^{p}}} \ \ a.s.
	\]
	and
	\[
	\lim_{t\rightarrow \infty} \langle x \rangle_{t} =  \dfrac{{(\delta_{1}^{l})^{1-p} (\delta_{1}^{u})^{p}}}{{(m_{1}^{l})^{1-p} (m_{1}^{u})^{p}}}( {(D^{l})^{1-p} (D^{u})^{p}} + \beta_{2})(R_{0}^{s} - 1) \ \ a.s.
	\]
\end{Theorem}

\noindent \textbf{Proof.}
According to (\ref{4.5}), we know that
\begin{equation}\label{4.10}
\begin{split}
\dfrac{\ln x(t)}{t}
\leq& {(m_{1}^{l})^{1-p} (m_{1}^{u})^{p}} S^{0} - {(D^{l})^{1-p} (D^{u})^{p}} - \beta_{2}- \dfrac{{(m_{1}^{l})^{1-p} (m_{1}^{u})^{p}}}{{(\delta_{1}^{l})^{1-p} (\delta_{1}^{u})^{p}}}\langle x \rangle_{t}\\
& + {(m_{1}^{l})^{1-p} (m_{1}^{u})^{p}}\phi(t)  + \dfrac{{(m_{2}^{l})^{1-p} (m_{2}^{u})^{p}}(t)}{t} + \dfrac{\widetilde{M}_{2}(t)}{t} + \dfrac{\ln x(0)}{t},
\end{split}
\end{equation}
which implies
\begin{equation*}
\begin{split}
\ln x(t)
\leq &[{(m_{1}^{l})^{1-p} (m_{1}^{u})^{p}} S^{0} - {(D^{l})^{1-p} (D^{u})^{p}} - \beta_{2}]t - \dfrac{{(m_{1}^{l})^{1-p} (m_{1}^{u})^{p}}}{{(\delta_{1}^{l})^{1-p} (\delta_{1}^{u})^{p}}}\int_{0}^{t}x(s)ds\\
& + {(m_{1}^{l})^{1-p} (m_{1}^{u})^{p}}\phi(t)t  + {(m_{2}^{l})^{1-p} (m_{2}^{u})^{p}}(t) + \widetilde{M}_{2}(t) + \ln x(0).
\end{split}
\end{equation*}
This together with (\ref{4.2}), (\ref{4.6}) and Zhao et al. (2013) implies
\begin{equation}\label{4.11}
\limsup_{t\rightarrow \infty} \langle x \rangle_{t} < \dfrac{{(\delta_{1}^{l})^{1-p} (\delta_{1}^{u})^{p}}}{{(m_{1}^{l})^{1-p} (m_{1}^{u})^{p}}} [ {(m_{1}^{l})^{1-p} (m_{1}^{u})^{p}} S^{0} - {(D^{l})^{1-p} (D^{u})^{p}} - \beta_{2} ] \ \ a.s.
\end{equation}
Thus, by (\ref{4.8}) we get
\begin{equation*}
\begin{split}
\limsup_{t\rightarrow \infty}\dfrac{\ln y(t)}{t} \leq & {(m_{2}^{l})^{1-p} (m_{2}^{u})^{p}}[\dfrac{{(\delta_{1}^{l})^{1-p} (\delta_{1}^{u})^{p}}}{{(m_{1}^{l})^{1-p} (m_{1}^{u})^{p}}} ( {(m_{1}^{l})^{1-p} (m_{1}^{u})^{p}} S^{0} \\
&- {(D^{l})^{1-p} (D^{u})^{p}} - \beta_{2} ) ]- {(D^{l})^{1-p} (D^{u})^{p}} - \beta_{3}\\
=& [\dfrac{{(m_{2}^{l})^{1-p} (m_{2}^{u})^{p}}{(\delta_{1}^{l})^{1-p} (\delta_{1}^{u})^{p}}}{{(m_{1}^{l})^{1-p} (m_{1}^{u})^{p}}}({(D^{l})^{1-p} (D^{u})^{p}} + \beta_{2}) \\
&+ {(D^{l})^{1-p} (D^{u})^{p}} + \beta_{3}] (R_{1}^{s} - 1)\\
<& 0 \ \ a.s. 
\end{split}
\end{equation*}
which yields
\begin{equation}\label{4.12}
\lim_{t\rightarrow \infty} y(t) = 0  \ \ a.s.
\end{equation}
It then follows from (\ref{4.5}) that
\begin{equation*}
\begin{split}
\langle x \rangle_{t} =& \dfrac{{(\delta_{1}^{l})^{1-p} (\delta_{1}^{u})^{p}}}{{(m_{1}^{l})^{1-p} (m_{1}^{u})^{p}}}[ {(m_{1}^{l})^{1-p} (m_{1}^{u})^{p}} S^{0} - {(D^{l})^{1-p} (D^{u})^{p}} - \beta_{2}  ]\\
& - ( \dfrac{1}{{(\delta_{2}^{l})^{1-p} (\delta_{2}^{u})^{p}}} + \dfrac{{(m_{2}^{l})^{1-p} (m_{2}^{u})^{p}}{(\delta_{1}^{l})^{1-p} (\delta_{1}^{u})^{p}}}{{(m_{1}^{l})^{1-p} (m_{1}^{u})^{p}}{(\delta_{2}^{l})^{1-p} (\delta_{2}^{u})^{p}}}) \langle y \rangle_{t}\\
& + \dfrac{{(\delta_{1}^{l})^{1-p} (\delta_{1}^{u})^{p}}}{{(m_{1}^{l})^{1-p} (m_{1}^{u})^{p}}} [  {(m_{1}^{l})^{1-p} (m_{1}^{u})^{p}}\phi(t)  + \dfrac{{(m_{2}^{l})^{1-p} (m_{2}^{u})^{p}}(t)}{t}\\
& + \dfrac{\widetilde{M}_{2}(t)}{t} + \dfrac{\ln x(0)}{t}    ]
- \dfrac{{(\delta_{1}^{l})^{1-p} (\delta_{1}^{u})^{p}}}{{(m_{1}^{l})^{1-p} (m_{1}^{u})^{p}}}\dfrac{\ln x(t)}{t}.
\end{split}
\end{equation*}
This together with Lemma~\ref{L42} implies
\begin{equation}\label{4.13}
\begin{split}
\liminf_{t\rightarrow \infty}\langle x \rangle_{t} =& \dfrac{{(\delta_{1}^{l})^{1-p} (\delta_{1}^{u})^{p}}}{{(m_{1}^{l})^{1-p} (m_{1}^{u})^{p}}} [ {(m_{1}^{l})^{1-p} (m_{1}^{u})^{p}}S^{0} - {(D^{l})^{1-p} (D^{u})^{p}} - \beta_{2}] \\
&- (\dfrac{1}{{(\delta_{2}^{l})^{1-p} (\delta_{2}^{u})^{p}}}+\dfrac{{(m_{2}^{l})^{1-p} (m_{2}^{u})^{p}}{(\delta_{1}^{l})^{1-p} (\delta_{1}^{u})^{p}}}{{(m_{1}^{l})^{1-p} (m_{1}^{u})^{p}}{(\delta_{2}^{l})^{1-p} (\delta_{2}^{u})^{p}}}) \limsup_{t\rightarrow \infty}\langle y \rangle_{t}\\
&+ \dfrac{{(\delta_{1}^{l})^{1-p} (\delta_{1}^{u})^{p}}}{{(m_{1}^{l})^{1-p} (m_{1}^{u})^{p}}}\liminf_{t\rightarrow \infty} [  {(m_{1}^{l})^{1-p} (m_{1}^{u})^{p}}\phi(t)  + \dfrac{{(m_{2}^{l})^{1-p} (m_{2}^{u})^{p}}(t)}{t}\\ 
&+ \dfrac{\widetilde{M}_{2}(t)}{t} + \dfrac{\ln x(0)}{t}    ]
-
\dfrac{{(\delta_{1}^{l})^{1-p} (\delta_{1}^{u})^{p}}}{{(m_{1}^{l})^{1-p} (m_{1}^{u})^{p}}}\limsup_{t\rightarrow \infty}\dfrac{\ln x(t)}{t}\\
\geq & \dfrac{{(\delta_{1}^{l})^{1-p} (\delta_{1}^{u})^{p}}}{{(m_{1}^{l})^{1-p} (m_{1}^{u})^{p}}}[ {(m_{1}^{l})^{1-p} (m_{1}^{u})^{p}}S^{0} - {(D^{l})^{1-p} (D^{u})^{p}} - \beta_{2} ] \ \ a.s.
\end{split}
\end{equation}
Therefore, from (\ref{4.11} ) and ( \ref{4.13}) we have
\begin{equation}\label{4.14}
\begin{split}
\lim_{t\rightarrow \infty}\langle x \rangle_{t} =& \dfrac{{(\delta_{1}^{l})^{1-p} (\delta_{1}^{u})^{p}}}{{(m_{1}^{l})^{1-p} (m_{1}^{u})^{p}}} [ {(m_{1}^{l})^{1-p} (m_{1}^{u})^{p}}S^{0} - {(D^{l})^{1-p} (D^{u})^{p}} - \beta_{2}]\\
= & \dfrac{{(\delta_{1}^{l})^{1-p} (\delta_{1}^{u})^{p}}}{{(m_{1}^{l})^{1-p} (m_{1}^{u})^{p}}} ({(D^{l})^{1-p} (D^{u})^{p}} + \beta_{2}) (R_{0}^{s} -1 ) \ \ a.s.
\end{split}
\end{equation}
It follows from (\ref{4.1}), (\ref{4.2}), (\ref{4.12}) and (\ref{4.14}) that 
\[
\lim_{t\rightarrow \infty} \langle S \rangle_{t} = \dfrac{ {(D^{l})^{1-p} (D^{u})^{p}} + \beta_{2}}{{(m_{1}^{l})^{1-p} (m_{1}^{u})^{p}}}\ \ a.s.
\]
This completes the proof.  \qed

\section{Persistence}
\label{sec:5}

In this section, we establish sufficient conditions for persistence in the mean of system (\ref{2.3}).
\begin{Definition}\label{D51}
	Model (\ref{2.3}) is said to be persistence in the mean, if
	\[
	\liminf_{t\rightarrow \infty}\langle y \rangle_{t} >0 \ \ a.s.
	\]
\end{Definition}

\begin{Theorem}\label{T51}
	Let Assumption (H3) hold. For any initial value $ (S(0), x(0), y(0))$ $\in \mathbb{R}^{3}_{+} $, there is a unique positive solution $ (S(t), x(t), y(t)) $ to system (\ref{2.3}) and the solution will remain in $\mathbb{R}^{3}_{+} $ with probability one, i.e., $ (S(t), x(t), y(t))\in \mathbb{R}^{3}_{+} $ for all $ t \geq 0 $ a.s. If $ R_{1}^{s}>1 $, then
	\begin{equation}\label{5.1}
	\begin{split}
	\liminf_{t\rightarrow \infty}\langle y \rangle_{t} \geq& \dfrac{{(m_{1}^{l})^{1-p} (m_{1}^{u})^{p}}{(\delta_{2}^{l})^{1-p} (\delta_{2}^{u})^{p}}}{{(m_{1}^{l})^{1-p} (m_{1}^{u})^{p}}{(m_{2}^{l})^{1-p} (m_{2}^{u})^{p}}+({(m_{2}^{l})^{1-p} (m_{2}^{u})^{p}})^{2}{(\delta_{1}^{l})^{1-p} (\delta_{1}^{u})^{p}}}\\
	&[\dfrac{{(m_{2}^{l})^{1-p} (m_{2}^{u})^{p}}{(\delta_{1}^{l})^{1-p} (\delta_{1}^{u})^{p}}}{{(m_{1}^{l})^{1-p} (m_{1}^{u})^{p}}}( {(D^{l})^{1-p} (D^{u})^{p}} + \beta_{2} )\\
	+& {(D^{l})^{1-p} (D^{u})^{p}} + \beta_{3}]
	(R_{1}^{s}-1)\\
	>&0 \ \ a.s.
	\end{split}
	\end{equation}
\end{Theorem}

\noindent \textbf{Proof.}
According to (\ref{4.8}), we know that
\begin{equation*}
\begin{split}
&\dfrac{{(m_{1}^{l})^{1-p} (m_{1}^{u})^{p}}}{{(m_{2}^{l})^{1-p} (m_{2}^{u})^{p}}{(\delta_{1}^{l})^{1-p} (\delta_{1}^{u})^{p}}}\dfrac{\ln y(t)}{t}\\ =&\dfrac{{(m_{1}^{l})^{1-p} (m_{1}^{u})^{p}}}{{(\delta_{1}^{l})^{1-p} (\delta_{1}^{u})^{p}}} \langle x \rangle_{t}
- \dfrac{{(m_{1}^{l})^{1-p} (m_{1}^{u})^{p}}}{{(m_{2}^{l})^{1-p} (m_{2}^{u})^{p}}{(\delta_{1}^{l})^{1-p} (\delta_{1}^{u})^{p}}}({(D^{l})^{1-p} (D^{u})^{p}}\\
& + \beta_{3} - \dfrac{M_{3}(t)}{t} - \dfrac{\widetilde{M}_{3}(t)}{t} - \dfrac{\ln y(0)}{t}).
\end{split}
\end{equation*}
This together with (\ref{4.5}) implies
\begin{equation*}
\begin{split}
&\dfrac{\ln x(t)}{t} + \dfrac{{(m_{1}^{l})^{1-p} (m_{1}^{u})^{p}}}{{(m_{2}^{l})^{1-p} (m_{2}^{u})^{p}}{(\delta_{1}^{l})^{1-p} (\delta_{1}^{u})^{p}}}\dfrac{\ln y(t)}{t} \\
=& {(m_{1}^{l})^{1-p} (m_{1}^{u})^{p}} S^{0} - {(D^{l})^{1-p} (D^{u})^{p}} - \beta_{2}\\
&- \dfrac{{(m_{1}^{l})^{1-p} (m_{1}^{u})^{p}}}{{(m_{2}^{l})^{1-p} (m_{2}^{u})^{p}}{(\delta_{1}^{l})^{1-p} (\delta_{1}^{u})^{p}}}({(D^{l})^{1-p} (D^{u})^{p}} + \beta_{3} )\\
&-(\dfrac{{(m_{1}^{l})^{1-p} (m_{1}^{u})^{p}}}{{(\delta_{1}^{l})^{1-p} (\delta_{1}^{u})^{p}}{(\delta_{2}^{l})^{1-p} (\delta_{2}^{u})^{p}}} + \dfrac{{(m_{2}^{l})^{1-p} (m_{2}^{u})^{p}}}{{(\delta_{2}^{l})^{1-p} (\delta_{2}^{u})^{p}}})\langle y \rangle_{t} + {(m_{1}^{l})^{1-p} (m_{1}^{u})^{p}}\phi(t)\\ 
&+\dfrac{{(m_{2}^{l})^{1-p} (m_{2}^{u})^{p}}(t)}{t} + \dfrac{\widetilde{M}_{2}(t)}{t} + \dfrac{\ln x(0)}{t}\\
& + \dfrac{{(m_{1}^{l})^{1-p} (m_{1}^{u})^{p}}}{{(m_{2}^{l})^{1-p} (m_{2}^{u})^{p}}{(\delta_{1}^{l})^{1-p} (\delta_{1}^{u})^{p}}} (\dfrac{M_{3}(t)}{t} + \dfrac{\widetilde{M}_{3}(t)}{t} + \dfrac{\ln y(0)}{t}),
\end{split}
\end{equation*}
which yields
\begin{equation}\label{5.2}
\begin{split}
\langle y \rangle_{t} =& \dfrac{{(\delta_{1}^{l})^{1-p} (\delta_{1}^{u})^{p}}{(\delta_{2}^{l})^{1-p} (\delta_{2}^{u})^{p}}}{{(m_{1}^{l})^{1-p} (m_{1}^{u})^{p}}+{(m_{2}^{l})^{1-p} (m_{2}^{u})^{p}}{(\delta_{1}^{l})^{1-p} (\delta_{1}^{u})^{p}}}[ {(m_{1}^{l})^{1-p} (m_{1}^{u})^{p}} S^{0}\\
& - {(D^{l})^{1-p} (D^{u})^{p}} - \beta_{2}-  \dfrac{{(m_{1}^{l})^{1-p} (m_{1}^{u})^{p}}}{{(m_{2}^{l})^{1-p} (m_{2}^{u})^{p}}{(\delta_{1}^{l})^{1-p} (\delta_{1}^{u})^{p}}}({(D^{l})^{1-p} (D^{u})^{p}} + \beta_{3} )  ]\\
& + \dfrac{{(\delta_{1}^{l})^{1-p} (\delta_{1}^{u})^{p}}{(\delta_{2}^{l})^{1-p} (\delta_{2}^{u})^{p}}}{{(m_{1}^{l})^{1-p} (m_{1}^{u})^{p}}+{(m_{2}^{l})^{1-p} (m_{2}^{u})^{p}}{(\delta_{1}^{l})^{1-p} (\delta_{1}^{u})^{p}}}[{(m_{1}^{l})^{1-p} (m_{1}^{u})^{p}}\phi(t) \\
&+ \dfrac{{(m_{2}^{l})^{1-p} (m_{2}^{u})^{p}}(t)}{t} + \dfrac{\widetilde{M}_{2}(t)}{t} + \dfrac{\ln x(0)}{t}\\
& +  \dfrac{{(m_{1}^{l})^{1-p} (m_{1}^{u})^{p}}}{{(m_{2}^{l})^{1-p} (m_{2}^{u})^{p}}{(\delta_{1}^{l})^{1-p} (\delta_{1}^{u})^{p}}} (\dfrac{M_{3}(t)}{t} + \dfrac{\widetilde{M}_{3}(t)}{t} + \dfrac{\ln y(0)}{t}) ]
\\
& - \dfrac{{(\delta_{1}^{l})^{1-p} (\delta_{1}^{u})^{p}}{(\delta_{2}^{l})^{1-p} (\delta_{2}^{u})^{p}}}{{(m_{1}^{l})^{1-p} (m_{1}^{u})^{p}}+{(m_{2}^{l})^{1-p} (m_{2}^{u})^{p}}{(\delta_{1}^{l})^{1-p} (\delta_{1}^{u})^{p}}} [ \dfrac{\ln x(t)}{t} \\
&+ \dfrac{{(m_{1}^{l})^{1-p} (m_{1}^{u})^{p}}}{{(m_{2}^{l})^{1-p} (m_{2}^{u})^{p}}{(\delta_{1}^{l})^{1-p} (\delta_{1}^{u})^{p}}}\dfrac{\ln y(t)}{t} ].\\
\end{split}
\end{equation}
Taking the inferior limit on both sides of (\ref{5.2}) and combining with
Lemma~\ref{L42}, from (\ref{4.2}) and (\ref{4.6}) we have
\begin{equation*}
\begin{split}
\liminf_{t\rightarrow \infty}\langle y \rangle_{t} \geq& \dfrac{{(\delta_{1}^{l})^{1-p} (\delta_{1}^{u})^{p}}{(\delta_{2}^{l})^{1-p} (\delta_{2}^{u})^{p}}}{{(m_{1}^{l})^{1-p} (m_{1}^{u})^{p}}+{(m_{2}^{l})^{1-p} (m_{2}^{u})^{p}}{(\delta_{1}^{l})^{1-p} (\delta_{1}^{u})^{p}}}[ {(m_{1}^{l})^{1-p} (m_{1}^{u})^{p}} S^{0}\\
& - {(D^{l})^{1-p} (D^{u})^{p}} - \beta_{2}-  \dfrac{{(m_{1}^{l})^{1-p} (m_{1}^{u})^{p}}}{{(m_{2}^{l})^{1-p} (m_{2}^{u})^{p}}{(\delta_{1}^{l})^{1-p} (\delta_{1}^{u})^{p}}}\\
&({(D^{l})^{1-p} (D^{u})^{p}} + \beta_{3} )  ]\\
= &  \dfrac{{(m_{1}^{l})^{1-p} (m_{1}^{u})^{p}}{(\delta_{2}^{l})^{1-p} (\delta_{2}^{u})^{p}}}{{(m_{1}^{l})^{1-p} (m_{1}^{u})^{p}}{(m_{2}^{l})^{1-p} (m_{2}^{u})^{p}}+({(m_{2}^{l})^{1-p} (m_{2}^{u})^{p}})^{2}{(\delta_{1}^{l})^{1-p} (\delta_{1}^{u})^{p}}}\\
&[ \dfrac{{(m_{2}^{l})^{1-p} (m_{2}^{u})^{p}}{(\delta_{1}^{l})^{1-p} (\delta_{1}^{u})^{p}}}{{(m_{1}^{l})^{1-p} (m_{1}^{u})^{p}}} ({(D^{l})^{1-p} (D^{u})^{p}} + \beta_{2})\\
&+ {(D^{l})^{1-p} (D^{u})^{p}} + \beta_{3}](R_{1}^{s} - 1) \ \ a.s.
\end{split}
\end{equation*}
Therefore, by the condition $ R_{1}^{s}>1 $, we can easily obtain (\ref{5.1}). \qed


\begin{thebibliography}{99}
	{
		
		\bibitem{1} H.L. Smith, P.E. Waltman, The Theory of the Chemostat: Dynamics of Microbial Competition, Cambridge University Press, Cambridge, 1995.
		\bibitem{2} L. Chen, X. Meng, J. Jiao, Biological Dynamics, Science Press, Beijing, 1993.
		
		
		\bibitem{11} R. May, Stability and Complexity in Model Ecosystems, Princeton University Press, NJ, 2001.
		\bibitem{12} L. Imhof, S. Walcher, Exclusion and persistence in deterministic and stochastic chemostat models, J. Differential Equations 217 (2005) 26–53.
		\bibitem{13} F. Campillo, M. Joannides, I. Larramendy-Valverde, Approximation of the Fokker–Planck equation of the stochastic chemostat, Math. Comput.
		Simulation 99 (2011) 37–53.
		\bibitem{14} C. Xu, S. Yuan, An analogue of break-even concentration in a simple stochastic chemostat model, Appl. Math. Lett. 48 (2015) 62–68.
		\bibitem{15} D. Zhao, S. Yuan, Critical result on the break-even concentration in a single-species stochastic chemostat model, J. Math. Anal. Appl. 434 (2016) 1336–1345.
		\bibitem{16} Q. Zhang, D. Jiang, Competitive exclusion in a stochastic chemostat model with Holling type II functional response, J. Math. Chem. 54 (2016) 777–791.
		\bibitem{17} L. Wang, D. Jiang, A note on the stationary distribution of the stochastic chemostat model with general response functions, Appl. Math. Lett.
		73 (2017) 22–28.
		\bibitem{18} S. Sun, Y. Sun, G. Zhang, X. Liu, Dynamical behavior of a stochastic two-species Monod competition chemostat model, Appl. Math. Comput. 298 (2017) 153–170.
		\bibitem{19} M. Sun, Q. Dong, J. Wu, Asymptotic behavior of a Lotka–Volterra food chain stochastic model in the chemostat, Stoch. Anal. Appl. 35 (2017)	1–17.
		
		
		\bibitem{20} D. Pal, G.S. Mahaptra, G.P. Samanta, Optimal harvesting of prey-predator ststem with interval biological parameters: a bioeconomic model, Math. Biosci. 241 (2013) 181–187.
		\bibitem{21} Q. Wang, Z. Liu, X. Zhang, R.A. Cheke, Incorporating prey refuge into a predator–prey system with imprecise parameter estimates, Comput. Appl.
		Math. 36 (2017) 1067–1084.
		\bibitem{22} S. Sharma, G.P. Samanta, Optimal harvesting of a two species competition model with imprecise biological parameters, Nonlinear Dynam. 77 (2014)
		1101–1119.
		\bibitem{222} D. Kiouach, Y. Sabbar, Ergodic Stationary Distribution of a Stochastic Hepatitis B Epidemic Model with Interval-Valued Parameters and Compensated Poisson Process, Comput. Math. Method. M. (2020), https://doi.org/10.1155/2020/9676501
		\bibitem{23} P. Panja, S.K. Mondal, J. Chattopadhyay, Dynamical study in fuzzy threshold dynamics of a cholera epidemic model, Fuzzy Inf. Eng. 9 (2017) 381–401.
		\bibitem{24} A. Das, M. Pal, A mathematical study of an imprecise SIR epidemic model with treatment control, J. Appl. Math. Comput. (2017) http://dx.doi.org/10.1007/s12190-017-1083-6.
		\bibitem{25} X.R. Mao, Stochastic Differential Equations and Applications, second ed., Horwood, Chichester, UK, 2007.	
		\bibitem{26} R. Situ, Theory of Stochastic Differential Equations with Jumps and Applications, Springer, 2005.	
		\bibitem{27}  M. Gao, D. Jiang, T. Hayat, A. Alsaedi, Threshold behavior of a stochastic Lotka–Volterra food chain chemostat model with jumps, Phys. A 523 (2019) 191–203.	
	}
\end{thebibliography}
\end{document}